%% file: art7.tex
\documentclass{amsart}
\usepackage{amssymb,xspace,textcomp}
\usepackage[small]{diagrams}
\usepackage[pdfstartview=FitH]{hyperref}
\input{command}\arrowsUsual\theorems

\def\Rep{\mathcal R}
\begin{document}
\input{title}
\input{abstract}
\maketitle
\input{introduction}
\input{section01} 
\input{section02} 
\input{section03}
\input{section04} 
\input{section05} 
\input{section06} 
\input{section07} 

\bibliography{../tex/fullbib}
\bibliographystyle{../tex/hamsplain}
\end{document}

%% file: command.tex
\def\mat#1{\ensuremath{#1}\xspace}

\def\cA{\mat{\mathbb{A}}}   

\def\cF{\mat{\mathbb{F}}}
\def\cN{\mat{\mathbb{N}}}   
\def\cQ{\mat{\mathbb{Q}}}   
\def\cC{\mat{\mathbb{C}}}   

\def\cZ{\mat{\mathbb{Z}}}   

\def\lE{\mat{\mathcal{E}}}

\def\lH{\mat{\mathcal{H}}}

\def\lL{\mat{\mathcal{L}}}
\def\lM{\mat{\mathcal{M}}}

\def\lO{\mat{\mathcal{O}}}
\def\lP{\mat{\mathcal{P}}}
\def\lR{\mat{\mathcal{R}}}

\def\lX{\mat{\mathcal{X}}}
\def\lY{\mat{\mathcal{Y}}}

\let\Phitemp\Phi \def\Phi{\mat{\Phitemp}}
\let\Psitemp\Psi \def\Psi{\mat{\Psitemp}}
\let\etatemp\eta \def\eta{\mat{\etatemp}}

\def\la{\mat{\lambda}}
\def\vi{\mat{\varphi}}
\let\mutemp\mu\def\mu{\mat{\mutemp}}
\let\nutemp\nu\def\nu{\mat{\nutemp}}
\let\pitemp\pi\def\pi{\mat{\pitemp}}
\def\si{\mat{\sigma}}

\def\al{\mat{\alpha}}
\def\be{\mat{\beta}}

\def\Ga{\mat{\Gamma}}

\def\De{\mat{\Delta}}

\def\hi{\mat{\chi}}

\def\ksi{\mat{\xi}}
\let\xitemp\xi \def\xi{\mat{\xitemp}}

\def\te{\mat{\theta}}

\def\mrm@#1{\mat{\mathrm{#1}}}

\def\g {\mat{\mathfrak{g}}}  

\def\gh{\mat{\mathfrak{h}}}
\def\gm{\mat{\mathfrak{m}}}

\def\DMO{\DeclareMathOperator}

\DMO{\Hom}{Hom}
\DMO{\lHom}{\lH\mathit{om}}
\DMO{\Ext}{Ext}
\DMO{\lExt}{\lE\mathit{xt}}
\DMO{\End}{End}
\DMO{\Aut}{Aut}
\DMO{\Fun}{Fun}
\DMO{\Tor}{Tor}
\DMO{\ext}{ext}

\DMO{\Ob}{Ob}
\DMO{\Mor}{Mor}
\DMO{\im}{im}
\DMO{\coim}{coim}
\DMO{\coker}{coker}
\DMO{\Arr}{Arr}

\DMO{\Id}{Id}
\DMO{\add}{add} 
\DMO{\ind}{ind} 
\DMO{\pro}{pro} 
\DMO{\Map}{Map} %
\DMO{\Iso}{Iso} %
\DMO{\Isom}{Isom}%
\DMO{\Ind}{Ind}

\DMO{\Presh}{Presh}
\def\Sch{\mat{\mathbf{Sch}}}

\DMO\coalg{Coalg}
\DMO{\Rep}{Rep}
\DMO{\Cor}{Cor}

\DMO{\Mod}{Mod}
\DMO{\rad}{rad}
\DMO{\soc}{soc}
\DMO{\ann}{ann}

\DMO{\Spec}{Spec}
\DMO{\spec}{Spec}
\DMO{\Proj}{Proj}
\DMO{\supp}{supp}
\DMO{\Coh}{Coh}
\DMO{\coh}{Coh}
\DMO{\Qcoh}{QCoh}
\DMO{\QCoh}{QCoh}
\DMO{\Pic}{Pic}
\DMO{\Div}{Div}
\DMO{\ch}{ch}
\DMO{\Hilb}{Hilb}
\DMO{\Fitt}{Fitt}
\DMO{\Quot}{Quot}
\def\Gm{\mat{{{\mathbb G}_{\mathrm m}}}}
\DMO{\Gras}{Gr}
\DMO{\Flag}{Flag}

\DMO{\cone}{cone}
\DMO{\Tw}{Tw}
\DMO{\rank}{rk}
\DMO{\rk}{rk}
\DMO{\codim}{codim}
\DMO{\cov}{cov}
\DMO{\sgn}{sgn}
\DMO{\td}{td}
\DMO{\GL}{GL}
\DMO{\SL}{SL}
\def\gl{\mat{\mathfrak{gl}}}

\DMO\Der{Der}
\DMO\der{Der}
\DMO\coder{Coder}
\DMO{\diag}{diag}
\DMO{\HMod}{HMod} 
\DMO{\ad}{ad}
\DMO*{\colim}{colim}
\DMO*{\hocolim}{hocolim}
\DMO*{\holim}{holim}
\DMO{\Ho}{Ho}

\DMO{\har}{char}
\DMO{\sk}{sk}
\DMO{\cosk}{cosk}
\DMO{\Gal}{Gal}
\DMO{\tr}{tr}
\DMO{\Tr}{Tr}
\DMO{\Sh}{Sh}
\DMO{\Is}{Is} 
\DMO{\Hol}{Hol} 
\DMO{\Lie}{Lie} 
\DMO{\Res}{Res} 
\DMO{\irr}{irr} %
\DMO{\Irr}{Irr} %
\DMO{\Exp}{Exp} %
\DMO{\Log}{Log} %
\DMO{\Pow}{Pow}
\DMO{\mult}{mult} %
\DMO{\height}{ht} %
\DMO{\wt}{wt}
\DMO{\Vect}{Vect}
\DMO{\moda}{mod}

\newcommand{\GIT}{/\!\!/}

\def\dd{\mat{\partial}}
\def\iso{\simeq}
\def\ts{\otimes}

\def\conjug#1{\overline{#1}}
\def\sb{\subset}

\def\iff{\mat{\Leftrightarrow}}

\def\xx{\times}

\def\ms{\backslash} 
\def\pser#1{[\![#1]\!]} 
\def\lpoly{\pser}

\def\inv{^{-1}}
\def\dual{^\vee}

\def\st{^{s}}
\def\sst{^{ss}}

\newcommand{\qbinom}[2]{\mat{\genfrac{[}{]}{0pt}{}{#1}{#2}}}  

\def\set#1{\mat{\{ #1\}}}
\def\sets#1#2{\mat{\{ #1 \mid #2\}}}

\def\arrowsUsual{
\newarrow{TeXto}----{->}
\newarrow{TeXinto}C---{->}
\newarrow{TeXonto}----{->>}
\newarrow{TeXdashto}{}{dash}{}{dash}{->}
\def\ar{\rightarrow}
\def\emb{\hookrightarrow}
\def\mto{\mapsto}
\def\arr{\rTeXto}
\def\embb{\rTeXinto}
\newarrow{Eq}=====
    }

\newif\ifukr\ukrfalse
\newif\ifrus\rusfalse
\newif\ifger\gerfalse

\def\theorems{
\newcounter{nthr} 
\numberwithin{nthr}{section}
\newtheorem{thr}[nthr]{Theorem}
\newtheorem{prp}[nthr]{Proposition}
\newtheorem{lmm}[nthr]{Lemma}
\newtheorem{crl}[nthr]{Corollary}
\newtheorem{clm}[nthr]{Claim}
\newtheorem{conj}[nthr]{Conjecture}
\theoremstyle{definition}
\newtheorem{dfn}[nthr]{Definition}
\newtheorem{rmr}[nthr]{Remark}
\newtheorem{exm}[nthr]{Example}}

\def\ub#1{\mat{\overline{#1}}}  
\def\lqq{\lq\lq}
\def\rqq{\rq\rq\xspace}
\def\eq{\equation}\def\endeq{\endequation}
\def\ie{i.e.\ }
\def\cf{cf.\ }

%% file: title.tex
\title[]{Fermionic forms and quiver varieties}%

\author{Sergey Mozgovoy}%

\address{
Department of Mathematics and Informatics,
University of Wuppertal,
42119 Wuppertal,
Germany.}

\email{mozgov@math.uni-wuppertal.de}%

\thanks{}%


%% file: abstract.tex
\begin{abstract}
We prove a formula relating the fermionic forms 
and the Poincar\'e polynomials of quiver varieties
associated to a finite quiver. Applied to quivers of type ADE, our
result implies a version of the fermionic Lusztig conjecture.
\end{abstract}

%% file: introduction.tex
\section{Introduction}
Let $(\Ga,I)$ be a finite quiver without loops, $C$ be the
corresponding generalized Cartan matrix, $\g=\g(\Ga)$ be the
corresponding Kac-Moody algebra and $\gh\sb \g$ be its Cartan
subalgebra, see \cite{Kac3}. Let $(\al_i)_{i\in I}$ be simple
roots of \g, $Q$ be the root lattice of \g, and $P$ be the weight
lattice of \g. For any root vector $\al\in Q$, we define
$(\al^i)_{i\in I}\in\cZ^I$ by $\al=\sum_{i\in I}\al^i\al_i$. We
will usually identify $Q$ with $\cZ^I$ in this way.
For any weight
vector $\nu\in P$, we define $(\nu_i)_{i\in I}\in\cZ^I$ by
$\nu_i:=(\nu,\al_i)$, where $(-,-)$ is a standard non-degenerate
symmetric bilinear form on $\gh\dual$ (see \cite[Section 2.1]{Kac3}).
It will be always clear from the context if the root
vector is considered as a weight vector so that there will
be no confusion with the above notation.
For any $\nu\in P$ and $\al\in Q$, we have
$(\nu,\al)=\sum_{i\in I}\nu_i\al^i$. We define the \cZ-valued
quadratic form $T$ (called the Tits form) on the root lattice $Q$
by $T(\al):=\frac 12 (\al,\al)$, $\al\in Q$.
For any $\nu\in P$ and $\al\in Q_+$, denote
$$[\nu,\al]=\prod_{i\in I}[\nu_i,\al^i],\qquad [\infty,\al]=\prod_{i\in I}[\infty,\al^i],$$
where the $q$-binomial coefficients $[n,m]$ (respectively, $[\infty,m]$),
for $n\in\cZ,\,m\in\cN$, are defined in Section \ref{sec:combi}.
Let \lP be the set of partitions. For any $\tau=(\tau^i)_{i\in
I}\in\lP^I$, we define $\tau_k:=(\tau^i_k)_{i\in I}\in\cN^I$,
$k\ge 1$ and $|\tau|:=\sum_{k\ge1}\tau_k$. We will usually
consider $\tau_k$ as elements of $Q_+$ (by identification of $Q$
with $\cZ^I$).

The fermionic forms $m(\nu,\la,q)\in\cQ(q)$ are usually attached
to $\nu=(\nu^{(k)}\in P_+)_{k\ge1}$ and $\la\in P_+$ (see
\cite[Section 4]{HKOTY1}, \cite[Section 4.2]{HKOTT1}). In this
paper we will assume that $\nu^{(k)}=0$ for $k>1$ and thus $\nu$
will be considered as an element of $P$. The relation of the
definition below to the usual definition of the fermionic forms
will be discussed in Section \ref{sec:combi}. For any
$\nu,\,\la\in P_+$, define the fermionic form by \eq
m(\nu,\la,q):=\sum_\tau
q^{-(\nu,\tau_1)}\prod_{k\ge1}q^{T(\tau_k)} \big[\nu-\sum_{l=1}^k
\tau_l,\tau_k-\tau_{k+1}\big],
\endeq
where the sum is taken over all partitions $\tau\in\lP^I$ such
that
\begin{enumerate}
    \item $\nu-\sum_{l=1}^k \tau_l\in P_+$, for any $k\ge 1$,
    \item $|\tau|=\nu-\la$.
\end{enumerate}
It is clear that $m(\nu,\la,q)$ can be nonzero only for
$\la\le\nu$ (\ie $\nu-\la\in Q_+$). For any $\nu,\la\in P$, define
the second fermionic form (see \cite[(4.16)]{HKOTY1}) by \eq
n(\nu,\la,q):=\sum_{|\tau|=\nu-\la}
q^{-(\nu,\tau_1)}\prod_{k\ge1}q^{T(\tau_k)} \big[\nu-\sum_{l=1}^k
\tau_l,\tau_k-\tau_{k+1}\big],
\endeq
where the sum is taken over all partitions $\tau\in\lP^I$ of
$\nu-\la$. The only difference with the definition of
$m(\nu,\la,q)$ is that now we do not impose any positivity
conditions on $\nu-\sum_{l=1}^k\tau_l$. Again, it is clear that
$n(\nu,\la,q)$ can be nonzero only for $\la\le\nu$.

It was conjectured by Lusztig \cite{Lus2} that if $\Ga$ is of type
ADE then the fermionic forms are related to the Poincar\'e
polynomials of Nakajima's quiver varieties associated to \Ga. Let
$\lM(\al,\nu)$ be the quiver variety (see \cite{Nak1,Nak2} and
Section \ref{sec moduli of repres}), where $\al\in Q_+$ and $\nu\in P_+$. It is
known that $\lM(\al,\nu)$ is a smooth symplectic manifold of
dimension $2d(\al,\nu)$, where $d(\al,\nu)=(\nu,\al)-T(\al)$.
Define $P(\lM(\al,\nu),q)=\sum_{i\ge0}\dim H_c^{i}(\lM(\al,\nu),\cC)q^{i/2}$.
In this paper we prove

\begin{thr}[Fermionic Lusztig conjecture]\label{intro:thr1}
If $\Gamma$ is of type ADE then, for any $\nu\in P_+$, we have
$$\sum_{\al\in Q_+}q^{-d(\al,\nu)}P(\lM(\al,\nu),q)e^{\nu-\al}%
=\sum_{\la\in P}n(\nu,\la,q\inv)\ch M(\la),$$%
where $M(\la)$ is the Verma module over $\g(\Ga)$ of highest
weight \la.
\end{thr}

It should be noted that the original fermionic conjecture of
Lusztig is slightly different from our formulation. It says that,
for any $\al\in Q_+$ and $\nu\in P_+$, one has
$$q^{-d(\al,\nu)}P(\lM(\al,\nu),q)%
=\sum_{\la\in P_+}m(\nu,\la,q\inv)\dim L(\la)_{\nu-\al},$$ %
where $L(\la)$ is the irreducible module of highest weight \la.
Equivalently, for any $\nu\in P_+$,
$$\sum_{\al\in Q_+}q^{-d(\al,\nu)}P(\lM(\al,\nu),q)e^{\nu-\al}%
=\sum_{\la\in P_+}m(\nu,\la,q\inv)\ch L(\la).$$%
As explained in \cite{Nak3}, this formula is a $q$-analog of the
Kirillov-Reshetikhin conjecture. In that conjecture there also
occur the fermionic forms $m(\nu,\la,q)$ (for $q=1$) and
characters $\ch L(\la)$ of irreducible modules. However, the known
proofs (see, e.g., \cite{Hernandez1,Nak3}) deal actually with a
slightly modified formula, which has the fermionic forms $n(\nu,\la,1)$
instead of $m(\nu,\la,1)$ and the Verma modules $M(\la)$
instead of the irreducible modules $L(\la)$. We have made
analogous modifications in the fermionic Lusztig conjecture.
The equality $m(\nu,\la,q)=n(\nu,\la,q)$ for $\nu,\la\in P_+$
was conjectured in \cite[Conjecture 4.3]{HKOTT1}. This problem is
still open, as far as we know. On the other hand, it follows
from the above theorem and
the existence of the Weyl group action on the cohomologies of
quiver varieties \cite{Lus4,Nak1} that
$$n(\nu,w\cdot\la,q)=(-1)^{l(w)}n(\nu,\la,q),$$
where $w$ is an element of the Weyl group and
$w\cdot\la=w(\la+\rho)-\rho$. This was also a part
of \cite[Conjecture 4.3]{HKOTT1}.

The above theorem will be an easy corollary of a more general
result, relating the fermionic forms with the Poincar\'e
polynomials of quiver varieties for arbitrary quivers \Ga. First,
let us recall the Hausel's formula (see \cite[Theorem
5]{Haus1} and Theorem \ref{thr:hausel form}), giving
an explicit expression for the generating function of
the Poincar\'e polynomials of quiver varieties.
For any $\nu\in P$, define $r(\nu,q)\in\cQ(q)\lpoly{x_i|i\in I}$ by
\equation\label{eq:r_nu}
r(\nu,q\inv):=\sum_{\tau\in\lP^I}x^{|\tau|}q^{-(\nu,\tau_1)}\prod_{k\ge1}q^{T(\tau_k)}
\big[\infty,\tau_k-\tau_{k+1}\big],
\endequation
where $x^\al=\prod_{i\in I}x_i^{\al^i}$ for $\al\in Q_+$.
The Hausel's formula says that, for any $\nu\in P_+$, we have
$$\sum_{\al\in
Q_+}q^{-d(\al,\nu)}P(\lM(\al,\nu),q)x^\al=\frac{r(\nu,q)}{r(0,q)}.$$
Define the generating function
$$
n(\nu,q)=\sum_{\la\in P}n(\nu,\la,q)x^{\nu-\la}=%
\sum_{\tau\in\lP^I}
x^{|\tau|}q^{-(\nu,\tau_1)}\prod_{k\ge1}q^{T(\tau_k)}
\big[\nu-\sum_{l=1}^k \tau_l,\tau_k-\tau_{k+1}\big].
$$
Our main combinatorial result (see Theorem \ref{thr:maincombi})
says that, for any $\nu\in P$, we have
\eq
n(\nu,q\inv)=r(\nu,q)r(0,q\inv).
\endeq
This is a purely combinatorial formula that can be easily verified
on computer. The starting point of the paper was an empiric
observation of this formula. Using it together with the Hausel's
formula, we get

\begin{thr}\label{intro:thr2}
For any $\nu\in P_+$, we have
$$\sum_{\al\in
Q_+}q^{-d(\al,\nu)}P(\lM(\al,\nu),q)x^\al=\frac{n(\nu,q\inv)}{r(0,q)r(0,q\inv)}.$$
\end{thr}
\noindent This result is used then in Section \ref{sec:ferm:Lusztig:conj} to prove the fermionic conjecture of Lusztig.
\medskip

The paper is organized as follows. In Section \ref{sec:combi} we
give some combinatorial prerequisites including \la-rings and the
definition of fermionic forms. In Section \ref{sec:poly-count} we
recall the basic results concerning the polynomial-count
varieties. In Section~\ref{GIT} we discuss the geometric invariant
theory over an arbitrary ring and prove the existence of good
quotients in the situation that occurs when constructing the
moduli spaces of quiver representations. In Section \ref{sec
moduli of repres} we recall the construction of moduli spaces of
stable quiver representations and, in particular, the construction
of quiver varieties. In Section \ref{sec:hausel formula} we give a
new proof of Hausel's formula based on Hua's formula
counting the absolutely indecomposable quiver representations over
finite fields. We also discuss there the proof of the first Kac
conjecture due to Hausel. In Section \ref{sec:main} we give a
proof of Theorem \ref{intro:thr2}. In Section
\ref{sec:ferm:Lusztig:conj} we use Theorem \ref{intro:thr2}
together with Hua's formula to prove the fermionic conjecture of
Lusztig.

%% file: section01.tex
\section{\texorpdfstring{$q$}{q}-binomial coefficients and
\texorpdfstring{\la}{lambda}-rings}\label{sec:combi}
For any expression $f$ depending on $q$, we define the conjugation
$\conjug{f(q)}:=f(q\inv)$, if it makes sense. Throughout the
paper, we will use the \la-ring structure on the ring of power
series over $\cQ(q)$ (see \cite[Appendix]{M2}). We endow the ring
$\cQ(q)\pser{x_1,\dots, x_r}$ with the structure of a \la-ring in
terms of Adams operations by
$$\psi_n(f(q,x_1,\dots,x_r))=f(q^n,x_1^n,\dots,x_r^n),
\qquad f\in\cQ(q)\pser{x_1,\dots,x_r}.$$
This \la-ring structure also endows $\cQ(q)\pser{x_1,\dots, x_r}$
with \la-operations and \si-ope\-ra\-tions.
Let $\gm$ be the maximal ideal of $\cQ(q)\pser{x_1,\dots, x_r}$.
Define the map $\Exp:\gm\ar 1+\gm$ by the formula
$$\Exp(f)=\sum_{k\ge0}\si_k(f)=\exp\Big(\sum_{k\ge 1}\frac 1k\psi_k(f)\Big),$$
where the map $\exp$ (as well as the map $\log$ used below) is
defined as in \cite[Ch.II \S 6]{BourLie13}.
For any $f,\,g\in \gm$, we have
$$\Exp(f+g)=\Exp(f)\Exp(g).$$
The map $\Exp$ has an
inverse $\Log:1+\gm\ar\gm$ which is given by the formula of Cadogan (see
\cite{Cadogan1,Getz1,M2})
$$\Log(f)=\sum_{k\ge1}\frac {\mu(k)}k\psi_k(\log(f)),$$
where $\mu$ is the M\"obius function.

For any $n\in\cZ$ and
$m\in\cN$, we define the $q$-binomial coefficients by
$$[n,m]=\qbinom {n+m}m:=\frac{\prod_{k=1}^m(1-q^{n+k})}{\prod_{k=1}^m(1-q^{k})}\in\cQ(q),$$
$$[\infty,m]:=\frac{1}{\prod_{k=1}^m(1-q^{k})}\in\cQ(q).$$
These functions can also be defined in the following way. Define
the $q$-Pochhammer symbols (\cf \cite{KC1})
\begin{multline*}
(x;q)_\infty=(1-x)^\infty_q:=\prod_{k\ge0}(1-q^kx)
=\prod_{k\ge0}\Exp(q^kx)\inv\\
=\Exp\Big(\sum_{k\ge0}q^kx\Big)\inv
=\Exp\left(\frac{x}{q-1}\right)\in \cQ(q)\pser x.
\end{multline*}
and
$$(x;q)_n=(1-x)^n_q:=\frac{(x;q)_\infty}{(q^nx;q)_\infty}=
\begin{cases}
\prod_{k=0}^{n-1}(1-q^kx),& n\ge0,\\
\prod_{k=n}^{-1}(1-q^kx)\inv,& n<0.\\
\end{cases}$$
Then
$$[n,m]=\frac{(q;q)_{m+n}}{(q;q)_n(q;q)_m}=\frac{(q^{n+1};q)_m}{(q;q)_m},\qquad
[\infty,m]=\frac1{(q;q)_m}.$$

\begin{thr}[Heine's formula, see {\cite[Theorem 13.1]{KC1}}]
We have
$$\sum_{k\ge0}\frac{(t;q)_k}{(q;q)_k}x^k
=\frac{(tx;q)_\infty}{(x;q)_\infty}=\Exp\left(\frac{1-t}{1-q}x\right)$$
in $\cQ(q)[t]\pser x$.
\end{thr}

In particular, taking $t=0$ or $t=q^{n+1}$, we obtain

\begin{lmm}
We have
$$\sum_{k\ge 0}[\infty,k]x^k=\Exp\left(\frac{x}{1-q}\right),$$
and, for any $n\in\cZ$,
$$\sum_{k\ge 0}[n,k]x^k=\Exp\left(\frac{1-q^{n+1}}{1-q}x\right).$$
\end{lmm}

%
%

For any $\nu=(\nu_i)_{i=1}^r\in\cZ^r$ and
$\al=(\al^i)_{i=1}^r\in\cN^r$, we define
$$[\nu,\al]:=\prod_{i=1}^r[\nu_i,\al^i].$$
The same notation is used for a weight vector
$\nu\in P$ and a root vector $\al\in Q_+$ (see the introduction).

We discuss now the relation between the definition of
fermionic forms in the introduction and the usual definition
\cite[(4.3)]{HKOTY1}, \cite[(4.5)]{HKOTT1}. Given a sequence
$\nu=(\nu_k\in P_+)_{k\ge1}$ 
with almost all elements being zero and $\la\in P_+$, one
defines loc.cit.
$$m(\nu,\la,q):=\sum_{m=(m_k)} q^{c(m)}\prod_{k\ge 1}[p_k(m),m_k],$$
where the sum runs over all $m=(m_k\in Q_+)_{k\ge 1}$ (with almost
all elements being zero) satisfying
\begin{enumerate}
    \item $p_k(m):=\sum_{l\ge 1}\min\set{k,l}(\nu_l-m_l)\in P_+$,
    for any $k\ge1$,
    \item $\sum_{k\ge1}k(\nu_k-m_k)=\la$
\end{enumerate}
and where
$$c(m):=\sum_{k,l\ge1}\min\set{k,l}\frac 12(m_k,m_l)
-\sum_{k,l\ge1}\min\set{k,l}(\nu_k,m_l).$$

If we define $\pi_k:=\sum_{l\ge k}\nu_l$ and
$\tau_k=\tau_k(m):=\sum_{l\ge k}m_l$ for $k\ge 1$ then it holds
\begin{enumerate}
    \item $p_k(m)=\sum_{l=1}^k (\pi_l-\tau_l)$,
    \item $\sum_{k\ge1}(\pi_k-\tau_k)=\la$,
    \item $c(m)=\frac
    12\sum_{k\ge1}(\tau_k,\tau_k)-\sum_{k\ge1}(\pi_k,\tau_k)$.
\end{enumerate}

This follows from the following easy fact

\begin{lmm}
Let $\pi,\tau\in\lP$ be partitions and let $n_k:=\pi_k-\pi_{k+1}$,
$m_k:=\tau_k-\tau_{k+1}$ for $k\ge1$. Then we have
$$\sum_{k,l\ge1}\min\set{k,l}n_k m_l=\sum_{k\ge1}\pi_k\tau_k.$$
\end{lmm}

We can now rewrite the definition of the fermionic form:
$$m(\nu,\la,q)=\sum_\tau \prod_{k\ge 1}q^{T(\tau_k)-(\pi_k,\tau_k)}
[\sum_{l=1}^k (\pi_l-\tau_l),\tau_k-\tau_{k+1}],$$%
where the sum runs over all partitions $\tau\in\lP^I$ (we identify
$\cZ^I$ with the root lattice~$Q$) satisfying
\begin{enumerate}
    \item $\sum_{l=1}^k (\pi_l-\tau_l)\in P_+$ for all $k\ge1$,
    \item $\sum_{k\ge1}(\pi_k-\tau_k)=\la$.
\end{enumerate}
The definition of the fermionic form given in
the introduction is precisely the last one.

%% file: section02.tex
\def\Fr{F}
\section{Polynomial-count varieties}\label{sec:poly-count}
By an algebraic variety over a field $k$ we will understand
a separated scheme of finite type over $k$.
All the algebraic varieties of this section will be assumed to be quasi-projective.
Given an algebraic variety $X$ over a field $k$ and a prime number
$l\ne\har k$, we define $H_c^i(X,\cQ_l)=H_c^i(X_{\ub k},\cQ_l)$
to be the $i$-th group of etale cohomologies with
compact support of $X_{\ub k}=X\ts_k\ub k$, where
$\ub k$ is a separable closure of $k$.

\begin{dfn}
Given an algebraic variety $X$ over a field $k$ and a prime
number $l\ne\har k$, we define the Poincar\'e polynomial
$P(X,q)\in\cZ[q^{1/2}]$ of $X$ by the formula
$$P(X,q)=\sum_{i\ge0}(-1)^i\dim H_c^i(X,\cQ_l)q^{i/2}.$$
\end{dfn}

\begin{dfn}
An algebraic variety $X$ over a finite field $k$
is called strictly polynomial-count if there exists a polynomial
$P\in\cQ[q]$ such that,
for any finite extension $\cF_q/k$, we have
$$P(q)=\#X(\cF_q).$$
$X$ is called polynomial-count if it is strictly
polynomial-count over some finite field extension of $k$.
In both cases the polynomial $P$ is called the counting polynomial of $X$.
\end{dfn}

\begin{exm}
The scheme given by $x^2=2$ in $\cA^1_{\cF_3}$ is polynomial-count
but not strictly polynomial-count.
\end{exm}

\begin{rmr}
Let $X$ be an algebraic variety over a finite field $\cF_q$.
If $X$ is strictly polynomial-count with a counting
polynomial $P(q)=\sum_{i\ge0}a_iq^i$ then its zeta-function
equals
$$Z(X,t)=\exp\Big(\sum_{n\ge1}\frac 1n\#X(\cF_{q^n})t^n\Big)
=\Exp(tP)(q)=\prod_{i\ge0}(1-tq^i)^{-a_i}.$$
As $Z(X,t)$ is a rational function, it follows that $P$ has integer
coefficients.
Moreover, $X$ is strictly polynomial-count
if and only if its zeta-function has zeros and poles
only of the form $q^{k}$, $k\in\cZ$.
Using the Grothendieck-Lefschetz formula one can express
the zeta-function \cite[1.5.4]{Del1}
$$Z(X,t)=\prod_{i\ge0}\det(1-t\Fr|H^i_c(X,\cQ_l))^{(-1)^{i+1}},$$
where $\Fr\in\Gal(\ub{\cF_q}/\cF_q)$ is a Frobenius element.
This implies that if the eigenvalues of Frobenius $\Fr$ on
$H^*_c(X,\cQ_l)$ are all of the form $q^k$, $k\in\cZ$, then $X$ is
strictly polynomial-count. If the eigenvalues of Frobenius $\Fr$
on $H^*_c(X,\cQ_l)$ are all of the form $\ksi q^k$, where $\ksi$
is a root of unity and $k\in\cZ$, then $X$ is polynomial-count
(\cf \cite[Definition 2.6]{LK1}). The last formula also implies,
that if the cohomologies of $X$ are concentrated in even degrees
and $X$ is polynomial-count then the coefficients of the counting
polynomial are nonnegative.
\end{rmr}

\begin{dfn}
Let $X$ be an algebraic variety over a finite field $\cF_q$
and let $l\ne\har\cF_q$ be a prime number.
$X$ is called $l$-pure \cite{CBB} if, for any $i\ge0$,
the conjugates of the eigenvalues of Frobenius $\Fr$ on $H^i_c(X,\cQ_l)$
have absolute value $q^{i/2}$ (the eigenvalues of Frobenius are
algebraic numbers \cite[Theorem 1]{Del2}).
\end{dfn}

\begin{lmm}[\cf \protect{\cite[Lemma A.1]{CBB}}]
Let $X$ be an algebraic variety over a finite field $\cF_q$
and let $l\ne\har \cF_q$ be a prime number.
Assume that $X$ is $l$-pure and is polynomial-count with a counting
polynomial $P$.
Then the cohomologies with compact support of $X$
are concentrated in even degrees and
the Poincar\'e polynomial of $X$ equals $P$ and has nonnegative coefficients.
\end{lmm}
\begin{proof}
We can suppose that $X$ is strictly polynomial-count.
As we have mentioned already, the zeta-function
$$Z(X,t)=\prod_{i\ge0}\det(1-t\Fr|H^i_c(X,\cQ_l))^{(-1)^{i+1}}$$
has zeros and poles of the form $q^k$, $k\in\cZ$.
Using the purity of $X$, we deduce that $H_c^i(X,\cQ_l)=0$ for odd $i$.
Moreover, the eigenvalues of $\Fr$ on even cohomologies are of
the form $q^k$, $k\in\cZ$. This implies that the eigenvalues of $\Fr$
on $H_c^{2i}(X,\cQ_l)$ equal $q^i$. It follows that
$$Z(X,t)=\prod_{i\ge0}(1-tq^i)^{-\dim H^{2i}_c(X,\cQ_l)}.$$
But we have seen that $Z(X,t)=\prod_{i\ge0}(1-tq^i)^{-a_i}$, where
$P(q)=\sum_{i\ge0}a_iq^i$. This implies that
$P(q)=\sum_{i\ge0}\dim H_c^{2i}(X,\cQ_l)q^i=P(X,q)$.
\end{proof}

Let $F$ be a number field. All places of $F$ that we will consider
are assumed to be finite (we identify them with the maximal ideals of
the ring of integers $\lO_F$).

\begin{rmr}
Let $F$ be a number field and let \lO be its ring of integers.
Given an algebraic variety $X$ over $F$,
there exists a scheme $X'$ over some localization $\lO_f$
of the ring $\lO$ such that
$X'\ts_{\lO_f} F\iso X$.
This allows us to define the schemes
$X'_{k(v)}=X'\ts_{\lO_f}k(v)$ over residue fields $k(v)$
for almost all places $v$ of $F$ (``almost all'' means here all except a finite set).
If $X''$ is a different scheme over some
localization $\lO_g$ such that
$X''\ts_{\lO_g} F\iso X$ then $X'$ and $X''$ are
isomorphic over some localization of $\lO$ and therefore
$X'_{k(v)}\iso X''_{k(v)}$ for almost all places $v$ of $F$.
We will assume that some scheme $X'$ as above is fixed and denote the
schemes $X'_{k(v)}$ just by $X_{k(v)}$.
For any field extension $k/k(v)$, we define $X(k)=X_{k(v)}(k)$.
Again, this is well-defined for almost all places
$v$ of $F$.
\end{rmr}

\begin{rmr}
Let $X$ be an algebraic variety over a number field $F$ and let
$l$ be a prime number.
Given an embedding $i:\ub F\emb\cC$, consider the
scheme $X_\cC=X\ts_F\cC$ and endow the set $X_\cC(\cC)$ with the
structure of a complex space.
By the Base Change Theorem \cite[Theorem 1.8.7]{FK1}, we have
$$H^i_c(X,\cQ_l)\iso H^i_c(X_\cC,\cQ_l).$$
By the Comparison Theorem \cite[Theorem 1.11.6]{FK1}, we have
$$H^i_c(X_\cC,\cQ_l)\iso H^i_c(X_\cC(\cC),\cQ_l),$$
where on the right hand side we consider the usual cohomologies with
compact support.
These formulas imply, in particular, that the Poincar\'e polynomial
of $X$ is independent of the prime number $l$.
For almost all places $v$ of $F$, we have
$$\dim H^i_c(X,\cQ_l)=\dim H^i_c(X_{k(v)},\cQ_l),\quad i\ge0.$$
Indeed, let $X$ be defined over a localization $\lO_f$, where
$\lO$ is a ring of integers of $F$. Denote by $\pi:X\ar\Spec\lO_f$
the corresponding structure map. For any $i\ge0$ the sheaf
$R^i\pi_!\cQ_l$ is constructible \cite[Theorem 1.12.15]{FK1}. This
implies \cite[Proposition 1.12.10]{FK1} that this sheaf is lisse
(\ie the corresponding $l$-adic projective system consists of
locally-constant sheaves) over some open subset of $\Spec\lO_f$.
It follows that the dimensions of the stalks at all points of this
open set coincide and this proves our statement. We deduce that,
for almost all places $v$ of $F$, the Poincar\'e polynomials of
$X_{k(v)}$ and $X$ coincide.
\end{rmr}

\begin{dfn}\label{dfn:purity etc. over number field}
An algebraic variety $X$ over a number field $F$
is called strictly polynomial-count if, for almost all places $v$
of $F$, the algebraic variety $X_{k(v)}$ over a finite field $k(v)$
is strictly polynomial-count. An algebraic variety $X$ over a number field $F$
is called polynomial-count if $X$ is strictly polynomial-count
over a finite field extension of $F$. Given a prime number $l$,
$X$ is called $l$-pure if, for almost
all places $v$ of $F$, the algebraic variety $X_{k(v)}$ is $l$-pure.
$X$ is called pure if it is $l$-pure, for all prime numbers $l$.
\end{dfn}

\begin{rmr}
If $X$ is a polynomial-count algebraic
variety over a number field $F$ then $X_{k(v)}$ is
polynomial-count for almost all places $v$ of $F$.
We do not know if the converse statement is true.
\end{rmr}

\begin{thr}[Katz {\cite[Theorem 6.1.4]{HRV}}]
Let $X$ be a polynomial-count algebraic variety over a number
field $F$. Then the algebraic varieties $X_{k(v)}$ have an
identical counting polynomial $P\in\cZ[q]$ for almost all places
$v$ of $F$. The virtual Hodge polynomial \cite{DanKhov} of $X$ can
be expressed as
$$E(X;u,v)=P(uv).$$
\end{thr}

\begin{prp}[see \cite{CBB}]\label{prp:purity criterion}
Let $X$ be a nonsingular algebraic variety over a number field $F$.
Assume that there exists a $\Gm$-action on $X$ with the properties
\begin{enumerate}
    \item for every $x\in X$ there exists the limit $\lim_{t\to0}tx$;
    \item $X^{\Gm}$ is projective.
\end{enumerate}
Then $X$ is $l$-pure for any prime number $l$.
\end{prp}

\begin{prp}[Nakajima {\cite{CBB}}]
\label{prp:equality in K0}
Let $X$ be a nonsingular algebraic variety over a number field $F$
and let $\pi:X\ar\cA^1$ be a smooth morphism over $F$.
Assume that there exists a $\Gm$-action on $X$ with the properties
\begin{enumerate}
  \item for every $x\in X$ there exists the limit $\lim_{t\to0}tx$;
    \item $\pi:X_F\ar\cA^1_{F}$ is $\Gm$-equivariant
    with respect to some $\Gm$-action on $\cA^1_F$ of positive weight.
\end{enumerate}
Then, for every $\la\in F$, the algebraic varieties $\pi\inv(0)$
and $\pi\inv(\la)$ represent the same element in $K_0(\Sch/\ub Q)$
(see e.g. \cite[Katz]{HRV}). In particular, if $\pi\inv(\la)$ is
polynomial-count then so is $\pi\inv(0)$ and they have the same
counting polynomial.
\end{prp}




%% file: section03.tex
\section{Geometric invariant theory}\label{GIT}
Let $R$ be a noetherian ring and let $G$ be an affine group
scheme over $R$.

\begin{dfn}\label{dfn:good Quot}
A morphism of finite type $\pi:X\ar Y$ of schemes over $R$, where $X$ is
endowed with an action of group $G$,
is called a good quotient if
\begin{enumerate}
    \item $\pi$ is $G$-invariant, affine, and surjective;
    \item the natural homomorphism $\lO_Y\ar (\pi_*\lO_X)^G$ is an isomorphism;
    \item for any closed $G$-stable subset $Z$ of $X$, $\pi(Z)$ is closed in $Y$;
    \item for any geometric point $R\ar k$, \pi induces a bijection between
  $Y(k)$ and the quotient of $X(k)$ by the equivalence relation
  $$x\sim y\iff \ub{Gx}\cap\ub Gy\ne\emptyset\text{ in }X\ts_R k.$$
\end{enumerate}
It is called a universal (respectively, a uniform) good quotient
if it stays a good quotient
after any base change (respectively, any flat base change) $Y'\ar Y$.
A good quotient is called a geometric quotient if all the orbits
of the geometric points of $X$ are closed.
\end{dfn}

\begin{rmr}
Sometimes the condition that \pi is affine is omitted (\cf
\cite[Def.~0.6]{GIT}). If $\pi:X\ar Y$ is a good quotient and
$Z_1,\dots,Z_n$ are $G$-stable closed subsets of $X$ then
$$\pi\big(\bigcap_{i=1}^n Z_i\big)=\bigcap_{i=1}^n\pi(Z_i).$$
It is proved in \cite[Remark 8]{S1} (\cf \cite[Remark 0.6]{GIT})
that good quotients are categorical quotients (\cite[Def.
0.5]{GIT}).
\end{rmr}

\begin{thr}[Seshadri \protect{\cite[Theorem 3]{S1}}]
\label{thr:seshadri}
Let $G$ be a reductive group scheme acting on an affine
scheme $X=\Spec A$ over $R$. Assume there exists a closed $G$-invariant
embedding $X\emb \cA^n_R=\Spec SV\dual$, where $V$ is a
$G$-$R$-module, free of rank $n$ over $R$. Then
$\pi:X\ar Y=\Spec A^G$ is a uniform good quotient.
If $R$ is universally Japanese
\cite[Ch.0, 23.1.1]{EGA4} then $Y$ is of finite
type over $R$.
\end{thr}
\begin{proof}
Everything except the last condition of Definition \ref{dfn:good Quot}
is proved in \cite[Theorem 3]{S1}.
The last condition is also stated there but it relies on
\cite[Proposition 6]{S1}, whose proof (and also the statement)
contains certain inaccuracies.
Still, we can prove a similar result. Namely, we prove
that for any two geometric points $x,y\in X(k)$ we have
$$\ub{Gx}\cap\ub{Gy}=\emptyset$$
in $X\ts_R k$ if and only if there exists some
$f\in\Ga(X,\lO_X)^G$ with $f(x)\ne f(y)$.
The \lqq if\rqq direction is trivial. To prove the
\lqq only if\rqq direction, we make the same reductions
as in Proposition \ref{prp:criterion semist} and assume that $X=\Spec SV\dual$,
$R$ is a discrete valuation ring or a field and $R\ar k$
is surjective.
Next, we can find $f_0\in SV\dual$ such that
$f_0(\ub{Gx})=1$ and $f_0(\ub{Gy})=0$.
The $R$-module $W$ generated by the $G$-translates
of $f_0$ in $S V\dual$ is finitely generated \cite[Prop. 3]{S1}
and is free as $R$ is a discrete valuation ring.
Consider a canonical map $\vi:\Spec SV\dual\ar\Spec SW$. Then
$\vi(x)$ is a $G$-invariant point in the affine space $\Spec(SW\ts_R k)$,
$\vi(y)=0$ and $\vi(x)\ne0$. Applying \cite[Theorem 1]{S1}, we can find
a homogeneous polynomial $f_1\in (SW)^G$ of positive degree such that
$f_1(\vi(x))\ne0$. It is clear that $f_1(\vi(y))=0$.
We define now $f=f_1\vi$. The result, which has just been proved,
implies the last
condition of Definition \ref{dfn:good Quot}
if we can show that $X(k)\ar Y(k)$ is
surjective. But $X$ is of finite type over $R$ and therefore
it is also of finite type over $Y$.
It follows that a surjective morphism $X\ar Y$ induces a
surjective map $X(k)\ar Y(k)$.
\end{proof}

\begin{rmr}
Excellent rings are universally Japanese \cite[7.8.3]{EGA4}.
In particular, any algebra of finite type over a field or
over \cZ is universally Japanese.
\end{rmr}

Let an affine group $G$ act on an affine scheme $X$ over $R$.
Given a character $\hi:G\ar \Gm$, we can endow the line bundle
$L=X\xx\cA^1_R$ over $X$ with a $G$-linearization, where the
action of $G$ on the second factor is given by \hi.
Let us recall \cite[Ch.1, Def. 1.7]{GIT} that a geometric
point $x$ of $X$ is called $L$-semistable (we will also say
$\hi$-semistable or just semistable) if there exists
a section $s\in\Ga(X,L^n)^G$, $n\ge1$ such that $s(x)\ne0$.
There exists an open subscheme $X\sst$ of $X$, such that the
semistable points of $X$ are precisely the geometric points
of $X\sst$. The embedding of graded algebras
$$\bigoplus_{n\ge0}\Ga(X,L^n)^G\ar\bigoplus_{n\ge0}\Ga(X,L^n)$$
induces the map
$$X\sst\ar\Proj\bigoplus_{n\ge0}\Ga(X,L^n)^G$$
which is known to be a good quotient if $R$ is an algebraically closed
field and $G$ is reductive \cite{King1}.
The goal of this section is to prove this result for an arbitrary ring $R$
(under certain conditions).
This result is not covered by the results of \cite{S1}
but, as we will see, it can be deduced from them. From now on,
we assume that $G$ is reductive and that there exists a closed $G$-invariant
embedding $X\emb \cA^n_R=\Spec SV\dual$, where $V$ is a
$G$-$R$-module, free of rank $n$ over $R$.

\begin{prp}[\cf \protect{\cite[Ch.0, 10.3.1]{EGA3}}, \protect{\cite[Theorem 1]{S1}}]
\label{prp:discrete extension}
Let $A$ be a discrete valuation ring, $m$ be its maximal
ideal, and $k=A/m$ be its residue field. Let $k\emb l$
be a field extension. Then there exists a flat homomorphism
$A\ar B$, where $B$ is a discrete valuation ring such that
the maximal ideal of $B$ equals $mB$ and
the residue field of $B$ is $k$-isomorphic to $l$.
\end{prp}
\begin{proof}
We will show that the construction of \cite[Ch.0, 10.3.1]{EGA3}
gives a discrete valuation ring if $A$ is a discrete valuation ring.

If $l=k(t)$, where $t$
is a variable, then one considers $A'=A[t]$, $m'=mA'$, and defines
$B=A'_{m'}$. It is proved in \cite[Ch.0, 10.3.1]{EGA3} that $B$ is local, noetherian
and satisfies all conditions of the proposition except the
condition that $B$ is a discrete valuation ring. But it is
clear that $B$ is a domain and its maximal ideal $m'B$ is a principal
ideal. This implies that $B$ is a discrete valuation ring.

Assume now that $l=k(t)$, where $t$ is algebraic over $k$. Let $f\in k[x]$
be the minimal polynomial of $t$ and let $F\in A[x]$ be the monic
polynomial that projects to~$f$. One defines then $B=A[x]/(F)$
and shows (see \cite[Ch.0, 10.3.1]{EGA3}) that $B$ is local, noetherian and
satisfies all conditions of the proposition except the fact that
$B$ is a discrete valuation ring. It is clear that the maximal ideal
$mB$ of $B$ is a principal ideal. To show that $B$ is a domain,
we note first that $F$ is irreducible in $K[x]$, where $K$ is a field
of fractions of $A$. This follows from Gau\ss's lemma
\cite[Ch.1, \S17, Lemma~1]{ZarSam1} and the fact that $f$ is irreducible
in $k[x]$. Consider now any two polynomials $F_1,F_2\in A[x]$
such that $F_1F_2\in FA[x]$. Then $F$ divides one of them in
$K[x]$ and therefore, again by Gau\ss's lemma, also in $A[x]$. This
implies that $A[x]/FA[x]$ is a domain.

The next steps of the proof repeat the remaining part of the proof
of \cite[Ch.0, 10.3.1]{EGA3}.
We just note that, given a filtered system $(A_i,f_{ij})$
of discrete valuation rings $(A_i,m_i)_{i\in I}$
such that $f_{ij}:A_i\ar A_j$ are flat homomorphisms and
$m_j=m_iA_j$ for $i<j$,
the colimit $B=\colim A_i$ of the system is a discrete valuation ring.
Indeed, it is a local noetherian ring with a maximal ideal
$m_B=m_iB$, for any $i\in I$, by \cite[Ch.0, 10.3.3]{EGA3}.
It follows that $m_B$ is a principal ideal. Moreover, $B$ is an
integral domain as a colimit of integral domains.
This implies that $B$ is a discrete valuation ring.
\end{proof}

The following result will be very important in Remark
\ref{rmr:king result}.

\begin{prp}\label{prp:criterion semist}
A geometric point $x\in X(k)$ is semistable if and only if
for some lift $\hat x\in L(k)$ of $x$, the closure $\ub{G\hat x}$
does not intersect the zero section of $L$.
\end{prp}
\begin{proof}
We can assume that $X=\Spec SV\dual$.
The ``only if'' part is trivial. Let us prove the ``if'' part.
Assume that there exists
some lift $\hat x\in L(k)$ of $x$ such that the closure $\ub{G\hat x}$
does not intersect the zero section.
Let $W$ be a $G$-module, free of rank $1$ over $R$, given by the character \hi,
and let $t\in W$ be its basis element.
Then $L\iso\Spec S(V\dual\oplus W\dual)$ as a $G$-scheme
and $\Ga(X,L^n)=SV\dual\ts W^{\ts n}$ as a $G$-module.
We have to show the existence of some section
$s=s_0\ts t^n\in (SV\dual\ts W^{\ts n})^G$, $n\ge1$ such that
$s(x)\ne0$ or, equivalently $s_0(x)\ne0$. Note that if $R\ar R'\ar k$ is a factorization
with a flat first morphism then it is enough to find the section
$s\in(SV\dual\ts W^{\ts n}\ts R')^{G\ts R'}$ with the same property
because
$$(SV\dual\ts W^{\ts n}\ts R')^{G\ts R'}=(SV\dual\ts W^{\ts n})^G\ts R'$$
by \cite[Lemma 2]{S1}.

After a flat base change, we may suppose that $G$ is split reductive
over $R$ (see \cite[Exp. 22, Cor. 2.3]{SGA3}, \cite[p.239]{S1}) and therefore $G$ is obtained
by a base change from a split reductive group $G_\cZ$ over $\cZ$
\cite[Exp. 25, Theorem 1.1]{SGA3}. By \cite[Prop. 1]{S1} we can embed
$V$ as a pure $G$-$R$-submodule into a finite direct sum of $R[G]$
(pure submodule means that it stays a submodule after any base change).
Then by \cite[Cor.2 of Prop.3]{S1} there exists a $G_\cZ$-module $U$, free
of finite rank over \cZ, and a pure embedding $V\emb U\ts_\cZ R$ of $G$-modules.
Hence, it is enough to find a section
$s\in (SU\dual\ts W_\cZ^{\ts n})^{G_\cZ}$ with $s(x)\ne0$. Thus, we may
suppose that $R=\cZ$.

If $\har k=0$ then $\cZ\ar k$ is flat and we may suppose by the above discussion
that $R=k$.
If $\har k=p>0$ then there exists a factorization $\cZ\ar \cZ_{(p)}\ar k$
with a flat first morphism. Moreover, by Proposition \ref{prp:discrete extension}
there exists a factorization $\cZ_{(p)}\ar B\ar k$ with a flat
first morphism and such that $B$ is a discrete valuation ring with a
residue field $k$. By the above discussion we can assume that $R=B$.
Thus, we can suppose that $R$ is a discrete valuation ring or a field
and $R\ar k$ is surjective.

There exists some $\ub f\in S(V\dual\oplus W)\ts k$ such that
$\ub f(\ub{G\hat x})=1$ and $\ub f(X\xx\set0)=0$.
The last equality means that
$\ub f\in SV\dual\ts(\bigoplus_{n\ge1}W^{\ts n})\ts k$.
Let us choose some representative
$f\in SV\dual\ts(\bigoplus_{n\ge1} W^{\ts n})$
of $\ub f$. The $R$-module $U$ generated by the $G$-translations
of $f$ in $SV\dual\ts(\bigoplus_{n\ge1}W^{\ts n})$ is finitely 
generated \cite[Prop. 3]{S1} and is free as $R$ is a discrete
valuation ring (or a field).
The induced morphism $\vi:L\ar \Spec SU$ maps $\hat x$ to a nonzero
$G$-invariant point. Applying \cite[Theorem 1]{S1},
we can find some $G$-invariant homogeneous polynomial $f_1\in SU$ of positive
degree such that $f_1(\vi(\hat x))\ne0$.
The image of $f_1$ under the map $SU\ar S(V\dual\oplus W)$ is contained in
$(SV\dual\ts(\bigoplus_{n\ge1}W^{\ts n}))^G$. Hence, for some
$s\in(SV\dual\ts W^{\ts n})^G$, $n\ge1$, we have $s(x)\ne0$.
\end{proof}

\begin{crl}\label{cor:sst under base change}
For any morphism $R\ar R'$, we have
$$(X\ts_R R')\sst=X\sst\ts_RR'.$$
\end{crl}

\begin{thr}\label{thr:semistable}
Let $R$ be a noetherian ring, $G$ be a reductive group scheme over $R$
acting on an affine scheme $X=\Spec A$ over $R$ and $\hi:G\ar\Gm$ be a
character.
Assume there exists a closed $G$-invariant
embedding $X\emb \cA^r_R=\Spec SV\dual$, where $V$ is a
$G$-$R$-module, free of rank $r$ over $R$. Then
the map
$$\pi:X\sst\ar Y=\Proj\bigoplus_{n\ge0}\Ga(X,L^n)^G,$$
where $X\sst$ is a subscheme of $\hi$-semistable points of $X$,
is a uniform good quotient.
If $R$ is universally Japanese
\cite[Ch.0, 23.1.1]{EGA4} then $Y$ is of finite
type over $R$ and the canonical morphism
$Y\ar\Spec A^G$ is projective.
\end{thr}
\begin{proof}
Let $B=\bigoplus_{n\ge0}\Ga(X,L^n)$. For any homogeneous
$f\in B^G$ of positive degree~$d$, consider an open subscheme
$Y_f=\Spec(B^G)_{(f)}$ of $Y$.
The preimage of $Y_f$ in $X\sst$ is $X_f=\Spec B_{(f)}$.
Note that $(B^G)_{(f)}=(B_{(f)})^G$ because
$(B^G)_f=(B_f)^G$ by \cite[Lemma 2]{S1}. Thus, it will
follow from Theorem \ref{thr:seshadri} that $X_f\ar Y_f$
is a good quotient if we show that there exists a closed
$G$-invariant embedding $X_f\emb\Spec SU$ for some $G$-module
$U$, free over $R$.
Let us denote by $W$ the $G$-module, free of rank~$1$ over $R$,
given by the character $\hi$. Then $L=\Spec(A\ts SW\dual)$
as a $G$-scheme and $\Ga(X,L^n)=A\ts W^{\ts n}$ as a $G$-module. It follows
$B\iso A\ts SW=A[t]$ as $G$-algebras,
where $t$ denotes some basis element of $W$.
The homogeneous element $f\in B^G$ can be written in the form $f_0t^d$,
where $f_0\in A$ satisfies $gf_0=\hi^{-d}(g)f_0$ for $g\in G$.
A surjective morphism of algebras
$A\ts S(W^{\ts d})=A[t^d]\ar B_{(f)}$, $t^d\mto 1/f_0$ is
$G$-invariant and defines the required embedding of $X_f$ in
$X\xx\cA^1\emb \cA^n\xx\cA^1$.

As $Y_f$ with homogeneous $f\in B^G$ cover $Y$, we obtain that
$X\sst\ar Y$ satisfies all properties (except the last one)
of a good quotient. To prove the last property of a good quotient,
we note first that $X\sst(k)\ar Y(k)$
is surjective, as we can apply Theorem \ref{thr:seshadri} locally.
If two geometric points in $X\sst(k)$ have different
images in $Y(k)$ then
the closures of their orbits do not intersect. Conversely,
assume that $x_1,x_2\in X\sst(k)$ are two geometric points
such that the closures of their orbits do not intersect.
We can find two homogeneous elements $f_1,f_2\in B^G$
of the same positive degree such that $f_1(x_1)\ne0$
and $f_2(x_2)\ne0$. Then one of the functions
$f_1,f_2,f_1+f_2$ is nonzero at both points $x_1,x_2$.
Let us denote this function by $f$. Then the closures
of the orbits of $x_1,x_2$ do not intersect in $X_f\ts_R k$.
Applying Theorem \ref{thr:seshadri}, we obtain that $x_1$ and
$x_2$ are mapped to different points under the map $X_f(k)\ar Y_f(k)$.

Assume now that $R$ is universally Japanese.
We have seen that $B\iso A\ts SW$ as a $G$-algebra
and therefore, by Theorem \ref{thr:seshadri}, $B^G$ is of
finite type over $R$. This implies that $B^G$ is of finite
type over $A^G$ and therefore $\Proj B^G\ar\Spec A^G$
is projective. As $A^G$ is of finite type over $R$ by
Theorem \ref{thr:seshadri}, we obtain that $\Proj B^G$ is of finite
type over $R$.
\end{proof}

\begin{rmr}\label{rmr:stable points}
A geometric point $x\in X(k)$ is called stable if it is
semistable, $Gx$ is closed in $X\sst\ts k$
and the dimension of $Gx$ equals $\dim G\ts k$.
There exists an open subscheme $Y\st\sb Y$ such that if $X\st=\pi\inv(Y\st)$
then the geometric points of $X\st$ are precisely the stable points
of $X$ (see \cite[Remark 9]{S1}).
The map $\pi:X\st\ar Y\st$ is a geometric quotient and it
induces a bijection $X\st(k)/G(k)\ar Y\st(k)$,
for any algebraically closed field $k$ (with a map $R\ar k$).
\end{rmr}

\begin{prp}
Let $G$ be a geometrically connected affine algebraic group
acting on a scheme $X$ over $R$ in such a way that the stabilizers of
the geometric points of $X$ are connected. Let $k$ be a finite field
(with a map $R\ar k$) and let $F\in\Gal(\ub k/k)$ be the Frobenius
element. Then the map
$$X(k)/G(k)\ar (X(\ub k)/G(\ub k))^F$$
is bijective.
\end{prp}
\begin{proof}
To prove the injectivity, assume that $x,y\in X(k)$
are mapped to the same point, that is,
there exists some $g\in G(\ub k)$ with $gx=y$.
This implies
$$gx=y=F(y)=F(gx)=F(g)F(x)=F(g)x.$$
Let $H\sb G(\ub k)$ be the stabilizer of $x$. Then $H$ is an $F$-invariant
connected affine algebraic group and $g\inv F(g)\in H$.
It follows from the Lang-Steinberg theorem \cite{Mueller1, Stein1} 
that there exists
some $h\in H$ such that $g\inv F(g)=hF(h\inv)$. Then
$F(gh)=gh$ and $(gh)x=gx=y$. This implies the injectivity.

To prove the surjectivity, consider a representative
$x\in X(\ub k)$ of some element from $(X(\ub k)/G(\ub k))^F$.
It satisfies $F(x)=gx$ for some $g\in G(\ub k)$. By the Lang-Steinberg
theorem there exists some $h\in G(\ub k)$ such that $g\inv=h\inv F(h)$.
This implies
$$F(hx)=F(h)F(x)=(hg\inv)(gx)=hx$$
and therefore $hx\in X(\ub k)^F=X(k)$.
\end{proof}

\begin{rmr}
Let $\lR_\al=\lR_\al(\Ga)$  be the space of representations
of a quiver \Ga and let $k$ be a finite field.
It is proved in \cite[Lemma 5.2.1]{KrR1} and \cite[Lemma 5.3.2]{KrR1}
that the map
$$\lR_\al(k)/\GL_\al(k)\ar(\lR_\al(\ub k)/\GL_\al(\ub k))^F$$
is bijective. The proof in \cite{KrR1} is completely elementary.
\end{rmr}

\begin{crl}\label{crl:points of stable quotient}
With the notation of Theorem \ref{thr:semistable},
assume that $G$ is geometrically connected and acts freely on geometric
points of $X\st$.
Then $\pi:X\st\ar Y\st$ induces a bijection
$X\st(k)/G(k)\ar Y\st(k)$,
for any finite field $k$ (with a map $R\ar k$).
\end{crl}
\begin{proof}
Let $\ub k$ be the algebraic closure of $k$ and
let $F\in\Gal(\ub k/k)$ be the Frobenius element. We know
that $X\st(\ub k)/G(\ub k)\ar Y\st(\ub k)$ is a bijection.
This implies that $Y\st(k)=Y\st(\ub k)^F$ can be identified with
$(X\st(\ub k)/G(\ub k))^F$ and the claim follows from the above
proposition.
\end{proof}



%% file: section04.tex
\section{Moduli spaces of quiver representations}\label{sec moduli of repres}
Let $(\Ga,I)$ be a finite quiver,
where \Ga is the set of arrows and $I$ is the set of vertices.
We will usually denote the quiver just by \Ga.
For any arrow $h\in \Ga$,
we denote by $h'$ and $h''$ its source and target, respectively.
For any field $k$, we denote by $k\Ga$ the path algebra of \Ga
over $k$. A $k\Ga$-representation will be considered as a pair
$(V,x)$, where $V=\oplus_{i\in I}V_i$
is an $I$-graded $k$-vector space and $x=(x_h)_{h\in\Ga}$
consists of homomorphisms $x_h:V_{h'}\ar V_{h''}$.

We denote by $\ub\Ga$ the double quiver of $\Ga$, obtained from it by
adjoining reverse arrows for all arrows in $\Ga$. For any $h\in
\ub\Ga$, we denote by $\ub h$ the opposite arrow contained in $\ub\Ga$.
For any $\al=(\al^i)_{i\in I}\in\cZ^I$  and
$\nu=(\nu_i)_{i\in I}\in\cZ^I$,
we define
$\nu\cdot\al:=\sum_{i\in I}\nu_i\al^i$.

Given an $I$-graded vector space $V=\oplus_{i\in I}V_i$, we
define its dimension to be $\dim V:=(\dim V_i)_{i\in I}\in\cN^I$.
Given two $I$-graded vector spaces $V$ and $W$, we
denote the vector space of $I$-graded morphisms between them by
$\Hom_I(V,W)$.

Let $V$ be an $I$-graded vector space of dimension \al.
Define the space of representations
$$\lR_\al(\Ga,k):=\bigoplus_{h\in\Ga}\Hom(V_{h'},V_{h''}).$$
We can identify $\lR_\al(\ub \Ga,k)$ with $\lR_\al(\Ga,k)\oplus \lR_\al(\Ga,k)\dual$.
There is an obvious action of the group
$$\GL_\al(k):=\prod_{i\in I}\GL_{\al^i}(k)$$
on $\lR_\al(\Ga,k)$ and on $\lR_\al(\ub \Ga,k)$ (this
action can be factored through \sloppy
$G_\al(k)=\GL_\al(k)\!/\Gm(k)$, where
$\Gm$ is considered as a diagonal subgroup in $\GL_\al$).

\begin{rmr}
For any commutative ring $R$, we can construct a free $R$-module
$\lR_\al(\Ga,R)$ in the same way as above.
We define then an affine space
$$\lR_\al(\Ga)_R=\Spec S(\lR_\al(\Ga,R)\dual)$$
over $R$. The scheme $\lR_\al(\Ga)=\lR_\al(\Ga)_\cZ$
is endowed with an action of the algebraic group $\GL_\al$ over \cZ.
\end{rmr}

Now we recall some facts about
the moduli spaces of semistable representations of quivers
from the paper of King \cite{King1}. Consider some
$\te\in\cZ^I$, called further stability.

\begin{dfn}\label{def:stability}
We define a slope function $\mu=\mu_\te:\cN^I\ms\set0\ar\cQ$ by the formula
$$\mu_\te(\al)=\frac{\te\cdot\al}{\height\al}.$$
Given a field $k$, we call a $k\Ga$-representation $V$ to be
\te-semistable (respectively, \te-stable) if, for any
nonzero proper subrepresentation $U\sb V$, we have
$\mu_\te(\dim U)\le\mu_\te(\dim V)$
(respectively, $\mu_\te(\dim U)<\mu_\te(\dim V)$).
A stable $k\Ga$-representation is called absolutely stable
if it stays stable after any field extension of $k$ (semistable
representations always stay semistable \cite[Lemma 4.2]{Rei1}).
We denote by $\lR_\al(\Ga,k)\sst$ the subset of
$\lR_\al(\Ga,k)$ consisting of semistable representations.
Denote by $\lR_\al(\Ga,k)\st$ the subset of
$\lR_\al(\Ga,k)$ consisting of absolutely stable representations.
\end{dfn}

\begin{rmr}
Stability \te is called \al-generic (or \al-coprime) if
$$\mu_\te(\be)\ne\mu_\te(\al),\qquad 0<\be<\al.$$
In this case all semistable points in $\lR_\al(\Ga,k)$ are automatically
absolutely stable.
\end{rmr}

\begin{rmr}[\cf \protect{\cite[section 2.2]{Rei4}}]\label{rmr:original stability}
The original definition of stability of representations \cite{King1}
is slightly different from Definition \ref{def:stability}.
There, given an element $\te'\in\cZ^I$, one calls a
representation $V$ of \Ga to be $\te'$-semistable if
$\te'\cdot\dim V=0$ and for any
subrepresentation $U\sb V$, we have $\te'\cdot\dim U\ge0$.
We can relate this definition to Definition \ref{def:stability}
as follows. Let $\al\in\cN^I\ms\set0$ and consider some
stability $\te\in\cZ^I$ as in Definition \ref{def:stability}.
Let $\rho\in\cZ^I$ be given by $\rho_i=1$, $i\in I$.
Let $a=\te\cdot\al$, $b=\rho\cdot\al=\height\al$ and
let $\te'=a\rho-b\te\in\cZ^I$.
Then $\te'\cdot\al=0$ and for any $\be\in\cN^I\ms{0}$, the condition
$\te'\cdot\be\ge0$ is equivalent to
$$\frac{\te\cdot\be}{\height\be}\le\frac ab=\frac{\te\cdot\al}{\height\al},$$
that is to $\mu_\te(\be)\le\mu_\te(\al)$. It follows that for representations
of dimension \al, the condition of
\te-semistability given in Definition \ref{def:stability} is
equivalent to the condition of $\te'$-semistability given at the
beginning of this remark.
\end{rmr}

\begin{rmr}\label{rmr:king result}
It was shown in \cite{King1} that the stability condition on
$\lR_\al(\Ga)$ can be interpreted in terms of geometric invariant
theory (see Section \ref{GIT}). Namely, for any stability \te, one
constructs a character $\hi:G_\al\ar\Gm$ by the rule
$$(g_i)_{i\in I}\mto\prod_{i\in I}(\det g_i)^{-\te'_i},$$
where $\te'$ was defined in Remark \ref{rmr:original stability}.
Then the geometric points of the scheme
$\lR_\al(\Ga)$ are $\hi$-semistable
(respectively, $\hi$-stable)
if and only if they are semistable (respectively, stable)
according to Definition \ref{def:stability}.
This result was proved in \cite{King1} for $R$ being an
algebraically closed field. But the proof in \cite{King1} uses
the definition of \hi-semistability (respectively, \hi-stability)
as in Proposition \ref{prp:criterion semist}
(respectively, Remark \ref{rmr:stable points} and Proposition \ref{prp:criterion semist})
and therefore works for arbitrary rings.
It follows that there exist open subschemes $\lR_\al(\Ga)\sst$
(respectively, $\lR_\al(\Ga)\st$) of $\lR_\al(\Ga)$
such that their geometric points are
precisely the semistable (respectively, stable) points of $\lR_\al(\Ga)$.
This subschemes are well behaved with respect to the base change
by Corollary \ref{cor:sst under base change}.
For any field $k$, we have
$$\lR_\al(\Ga)\sst(k)=\lR_\al(\Ga,k)\sst,\qquad
\lR_\al(\Ga)\st(k)=\lR_\al(\Ga,k)\st.$$
We can use the results of the previous section
to construct the quotient of $\lR_\al(\Ga)\sst$ (or $\lR_\al(\Ga)\st$)
by the action of $G_\al$ over an arbitrary ring $R$.
\end{rmr}


\begin{dfn}
A $k\Ga$-representation $(V,x)$ is called nilpotent if there
exists some $N\ge1$ such that for any path $h_1\dots h_N$ in \Ga
(\ie $h''_i=h'_{i+1}$, $1\le i<N$), we have $x_{h_N}\dots
x_{h_1}=0$. Denote by $\lR_\al(\Ga,k)^{n}$ the subset of
$\lR_\al(\Ga,k)$ consisting of nilpotent representations. Define
also
$\lR_\al(\Ga,k)^{ns}=\lR_\al(\Ga,k)^{n}\cap\lR_\al(\Ga,k)^{s}$.
\end{dfn}

The action of $\GL_\al$ on
$\lR_\al(\ub\Ga)=\lR_\al(\Ga)\oplus\lR_\al(\Ga)\dual$
induces a moment map
$$\mu:\lR_\al(\ub \Ga)\ar \g_\al\dual\emb\gl_\al\dual$$
given by
$$(x_h)_{h\in\ub\Ga}\mto \sum_{h\in \Ga}[x_h,x_{\ub h}],$$
where $\gl_\al=\prod_i \gl_{\al^i}$ is a Lie algebra of $\GL_\al$,
$\g_\al$ is a Lie algebra of $G_\al$,
$\gl_\al\dual$ is identified with $\gl_\al$ by means
of the trace pairing and
$\g_\al\dual\emb\gl_\al\dual$ is identified with matrices
$(\ksi_i)_{i\in I}$ such that $\sum_i\tr\ksi_i=0$.

\begin{lmm}\label{lmm:smoothness of moment map}
The moment map
$\mu:\lR_\al(\ub\Ga)\st\ar\g_\al\dual$ is smooth.
\end{lmm}
\begin{proof}
Assume first that $R$ is an algebraically closed field.
Stability condition implies that the stabilizer in
$G_\al$ of any stable point is trivial. This implies that the differential
of $\mu$ is surjective on tangent spaces and therefore $\mu$ is smooth
at stable points.
The smoothness of \mu for an arbitrary field $R$ follows now from
a faithfully flat descent \cite[17.7.3]{EGA4}.
Applying \cite[Exp.2, Cor.2.2]{SGA1}, we obtain that \mu
is smooth for an arbitrary noetherian ring $R$.
\end{proof}

\begin{crl}
If $\mu\inv(0)\st$ is nonempty then the map
$\mu:\lR_\al(\ub \Ga)\st\ar\g_\al\dual$ is surjective.
\end{crl}
\begin{proof}
The smooth morphisms are open, so
the image of the map $\mu:\lR_\al(\ub \Ga)\st\ar\g_\al\dual$ is
open. As this image contains $0$ and is stable with respect to the
multiplication by scalars, it coincides with the whole $\g_\al\dual$.
\end{proof}

\begin{lmm}\label{lem equality for fibers}
Assume that the stability \te is \al-generic
and that there exists some $G_\al$-invariant element
$\ksi\in\g_\al\dual(\cQ)$ with
a nonempty $\mu\inv(\ksi)\st\sb\lR_\al(\Ga)_\cQ\st$. Then the map
$\mu:\lR_\al(\ub\Ga)_\cQ\st\ar\g_{\al,\cQ}\dual$ is surjective
and the varieties $\mu\inv(\ksi)\st\GIT G_\al$ and
$\mu\inv(0)\st\GIT G_\al$ represent the same element in
$K_0(\Sch/\ub\cQ)$.
\end{lmm}
\begin{proof}
In view of the last corollary,
we have to prove only the second statement .
Let $L$ be a line through points $0$ and $\ksi$ in the affine space $\g_{\al,\cQ}\dual$.
Consider the schemes $X=\mu\inv(L)\st$, $Y=\mu\inv(L)$, $\lX=X\GIT G_\al$ and $\lY:=Y\GIT G_\al$ over \cQ.
There is a commutative diagram \begin{diagram}
X&\rInto&Y\\
\dTo&&\dTo \\
\lX&\rTo^p&\lY&\rTo&L,
\end{diagram}
where $p$ is projective by Theorem \ref{thr:semistable}. The
composition $\lX\ar\lY\ar L$ is smooth. Indeed, \lX is smooth by
Luna's etale slice theorem and the differential of $\lX\ar L$ is
surjective on tangent spaces as this is so for $X\ar\lX\ar L$ by
Lemma \ref{lmm:smoothness of moment map}. Let $\Gm$ act on
$\lR_\al(\Ga)$ by multiplication of all the matrices by scalar.
This action induces the action of $\Gm$ on $X$, $Y$, $\lX$ and
$\lY$. There is also an action of $\Gm$ on $L$ such that the
rightmost map in the above diagram is $\Gm$-equivariant. The
assertion of the lemma will follow from Propostion
\ref{prp:equality in K0} if we show that for any $x\in \lX$
there exists $\lim_{t\to 0}tx$. The map $\Gm\ar \lY$, $t\mto
tp(x)$ can be extended to $\cA^1$ by $0\mto 0$. Now it follows
from the projectivity of $p$ that the map $\Gm\ar\lX$, $t\mto tx$
can also be extended to $\cA^1\ar\lX$. This proves the existence
of the limit.
\end{proof}

\begin{rmr}
Multiplying the element $\ksi$ by some integer, we may assume that
$\ksi\in\g_\al\dual(\cZ)$. Then all constructions of the above
lemma can be made over \cZ. Note that the formation of quotients
commutes with flat base change $\cZ\ar\cQ$.
\end{rmr}

\subsection{Quiver varieties}
Let $K$ be an algebraically closed field. For any $\nu\in P_+$,
define an enlarged quiver $\Ga_*$ by adjoining to $\Ga$ a new
vertex $*$ and $\nu_i$ arrows from $*$ to $i$, for each $i\in I$.
Its set of vertices is $I_*=I\cup\set{*}$. For any $\al\in\cZ^I$
and any $k\in\cZ$, we consider the pair $(\al,k)$ as an element of
$\cZ^{I_*}$.

Let $W$ be an $I$-graded vector space of dimension $(\nu_i)_{i\in
I}$ and let $V$ be an $I$-graded vector space of dimension \al. We
identify $W_i$ with $\bigoplus_{h:*\ar i}K h$, for $i\in I$.
Define
$$M(\al,\nu):=\lR_{(\al,1)}(\ub \Ga_*)=\lR_\al(\ub\Ga)\oplus\Hom_I(W,V)\oplus\Hom_I(V,W).$$
The elements of this space will be represented as triples
$(x,p,q)$ with elements corresponding to the above decomposition.
Note that $G_{(\al,1)}=(\prod_{i\in
I}\GL_{\al_i}\xx\Gm)/\Gm\iso\GL_\al$. Therefore the moment map can
be considered as a map $\mu_*:M(\al,\nu)\ar\gl_\al\dual$. It is
given by the formula
$$\mu_*(x,p,q)=\mu(x)+pq,$$
where we again identify $\gl_\al\dual$ with $\gl_\al$ by the trace
pairing.

Let us fix a stability $\te=(0,\dots,0,-1)\in\cZ^{I_*}$.

\begin{lmm}\label{lmm:stab of quiver var}
Stability \te is $(\al,1)$-generic. An element $(x,p,q)\in
M(\al,\nu)$ is stable if and only if any $I$-graded $x$-invariant
subspace $V'\sb V$ such that $q(V')=0$, is zero.
\end{lmm}
\begin{proof}
Any element $0<(\be,k)<(\al,1)$ is either of the form $(\be,0)$
with $0<\be\le\al$ or of the form $(\be,1)$ with $0\le\be<\al$. In
the first case we have $\mu_\te(\be,0)=0>\mu_\te(\al,1)$. In the
second case we have
$$\mu_\te(\be,1)=\frac{-1}{\height\be}<\frac{-1}{\height\al}=\mu_\te(\al,1).$$
This shows that $\te$ is $(\al,1)$-generic and that a stable
representation $(x,p,q)$ of dimension $(\al,1)$ does not contain
nonzero subrepresentations of dimension $(\be,0)$. These
representations correspond to $I$-graded $x$-invariant subspaces
$V'\sb V$ such that $q(V')=0$.
%
\end{proof}

\begin{dfn}
Define the quiver variety $\lM(\al,\nu)$ to be the good quotient
$\mu_*\inv(0)\st\GIT\GL_\al$. Define
$\lL(\al,\nu)=\mu_*\inv(0)^{ns}\GIT\GL_\al$.
\end{dfn}

\begin{rmr}
It is easy to see that $\lL(\al,\nu)$ is the preimage of zero
under the projective morphism $\mu_*\inv(0)\st\GIT \GL_\al\ar
\mu_*\inv(0)\GIT \GL_\al$. It is known that an element $(x,p,q)\in
M(\al,\nu)\st$ is nilpotent if and only if $x$ is nilpotent and
$p=0$, see e.g.\ \cite[Lemma 5.9]{Nak1} or \cite[Lemma
2.22]{Lus1}.
\end{rmr}

Let $T$ denote the Tits form of the quiver \Ga and let $T_*$ denote
the Tits form of the quiver $\Ga_*$. As in the introduction, we define
$d(\al,\nu):=1-T_*(\al,1)=\al\cdot\nu-T(\al)$.

\begin{thr}[{Nakajima \cite[Section 3]{Nak2}}]
The variety $\lM(\al,\nu)$ is smooth and the variety $\lL(\al,\nu)$ is
projective. The complex manifold $\lM(\al,\nu)(\cC)$ is symplectic
and its subvariety $\lL(\al,\nu)(\cC)$ is a Lagrangian subvariety
homotopic to $\lM(\al,\nu)(\cC)$. The dimension of $\lM(\al,\nu)$
equals $2d(\al,\nu)$.
\end{thr}

%% file: section05.tex
\section{Hausel Formula}\label{sec:hausel formula}
In this section we give a new proof of Hausel's formula based
on Hua's formula for the absolutely indecomposable representations
of a quiver.

Let $(\Ga,I)$ be a finite quiver.
It was shown by Kac \cite{Kac1} that, for any $\al\in\cN^I$, there
exist polynomials $a_\al(\Ga)\in\cZ[q]$  and $m_\al(\Ga)\in\cZ[q]$
such that for any finite field
$\cF_q$, $a_\al(\Ga,q)$ (respectively, $m_\al(\Ga,q)$)
is the number of isomorphism classes of absolutely
indecomposable representations (respectively, all representations)
of $\Ga$ over $\cF_q$ of dimension \al.
It was proved by Kac \cite{Kac1} that
$a_\al\ne 0$ if and only if $\al$ is a root of
$\g(\Ga)$ (the Kac-Moody algebra associated to \Ga)
and $a_\al=1$ if and only if $\al$ is a real root.
Let us define the generating functions
$a(\Ga,q)=\sum_{\al\in\cN^I}a_\al(\Ga,q)x^\al$
and $m(\Ga,q)=\sum_{\al\in\cN^I}m_\al(\Ga,q)x^\al$.
Then we have \cite[Lemma 5]{M2}
$$m(\Ga,q)=\Exp(a(\Ga,q)).$$
The formula of Hua (see \cite[Theorem 6]{M2}) says that
$$\Exp\left(\frac{a(\Ga,q)}{q-1}\right)=r(\Ga,q),$$
where $r(\Ga,q)=\sum_{\al\in\cN^I}r_\al(\Ga,q)x^\al$ and
$$r_\al(\Ga,q\inv)=\sum_{|\tau|=\al}\prod_{k\ge1}q^{T(\tau_k)}[\infty,\tau_k-\tau_{k-1}].$$

\begin{prp}[see \protect{\cite[Proposition 2.2.1]{CBB}}]
\label{prp number of points of qv} The quiver variety
$\lM=\lM(\al,\nu)$ is polynomial-count with a counting polynomial
$q^{d(\al,\nu)}a_{(\al,1)}(\Ga_*,q)$.
\end{prp}
\begin{proof}
Recall that our stability is $\te=(0,\dots,0,-1)$.
Applying the construction of Remark \ref{rmr:original stability},
we get an element $\te'=(-1,\dots,-1,\sum\al_i)$. We consider
it as an element of $\gl_{(\al,1)}\iso\gl_{(\al,1)}\dual$
consisting of diagonal matrices. Consider the moment map
$\mu_*:\lR_{(\al,1)}(\ub\Ga_*)\ar\gl_{(\al,1)}$.
All geometric points of $\mu_*\inv(\te')$ describe irreducible
(and in particular stable) representations. Indeed, let
$(x,p,q)$ be any such point and let $V'\oplus W'\sb V\oplus K$
be a nonzero proper $I_*$-graded subspace of dimension $(\be,k)$,
stable under $(x,p,q)$. Then $\mu_*((x,p,q)|_{V'\oplus W'})=\te'$
as an element of $\gl{(\be,k)}$. Hence
$$\te'\cdot(\be,k)=\te'\cdot(\al,1)=0,$$
as the sum of the traces of any element from the image of the moment map is zero.
This implies that $\mu_\te(\be,k)=\mu_\te(\al,1)$ which is impossible
by Lemma \ref{lmm:stab of quiver var}. Hence
$\mu_*\inv(\te')\st=\mu_*\inv(\te')$.

It was proved in \cite[Proposition 2.2.1]{CBB} that
$$\#(\mu_*\inv(\te')(\cF_q)/\GL_{(\al,1)}(\cF_q))=q^{d(\al,\nu)}a_{(\al,\nu)}(\Ga_*,q).$$
It follows from Corollary \ref{crl:points of stable quotient} that
the geometric quotient $\mu_*\inv(\te')\st\GIT\GL_{(\al,1)}$ is
polynomial-count with a counting polynomial as above.
Now Lemma \ref{lem equality for fibers} implies that
also $\mu_*\inv(0)\st\GIT\GL_{(\al,1)}$ is polynomial-count
with the same counting polynomial.
\end{proof}

We denote the counting polynomial of $\lM=\lM(\al,\nu)$ by
$P(\lM,q)$.

\begin{prp}
The quiver variety $\lM(\al,\nu)$ is pure (see Def. \ref{dfn:purity etc. over number field}).
Its Poincar\'e polynomial (of cohomologies with compact support) equals
$P(\lM,q)$.
\end{prp}
\begin{proof}
We have to check the conditions of Proposition
\ref{prp:purity criterion} for the quiver variety
$\lM=\mu_*\inv(0)\st\GIT\GL_\al$.
We define the same action of $\Gm$ on $\mu_*\inv(0)\st\GIT\GL_\al$
as in Lemma \ref{lem equality for fibers}.
The existence of the limit is proved in the same way as
there.
The \Gm-invariant part of $\lM$ is mapped to
zero under the projective morphism $p:\mu_*\inv(0)\st\GIT\GL_\al\ar\mu_*\inv(0)\GIT\GL_\al$,
as the only $\Gm$-invariant point of $\lR_{(\al,1)}\GIT\GL_\al$
is zero.
This implies that $\lM^\Gm$ is projective.
\end{proof}

\begin{thr}[Hausel formula]\label{thr:hausel form}
We have
$$\sum_{\al\in Q_+}q^{-d(\al,\nu)}P(\lM(\al,\nu),q)x^\al=\frac{r(\nu,q)}{r(0,q)},$$
where $r(\nu,q)$ is given by equation \eqref{eq:r_nu}.
\end{thr}
\begin{proof}
For any $n\in\cN$, we define
$$a_n:=\sum_{\al\in\cN^I}a_{(\al,n)}(\Ga_*)x^\al,\qquad %
r_n:=\sum_{\al\in\cN^I}r_{(\al,n)}(\Ga_*)x^\al,$$
$a_*:=\sum_{n\ge0}a_nx_*^n$ and $r_*:=\sum_{n\ge0}r_nx_*^n$.
Applying Hua's formula to the quiver $\Ga_*$ we get
$a_*(q)=(q-1)\Log(r_*(q))$ and therefore
\begin{multline*} a_1(q)=\frac \dd{\dd x_*}a_*(q)\big|_{x_*=0}
=(q-1)\frac \dd{\dd
x_*}\sum_{k\ge1}\frac {\mu(k)}k\psi_k(\log(r_*(q)))\big|_{x_*=0}\\
=(q-1)\frac \dd{\dd x_*}\log(r_*(q))\big|_{x_*=0}
=(q-1)\frac{\frac \dd{\dd x_*}r_*(q)}{r_*(q)}\big|_{x_*=0}
=(q-1)\frac {r_1(q)}{r_0(q)}.
\end{multline*}
We will show in the next lemma that $r(\nu,q)=(q-1)r_1(q)$. As
$r_0(q)=r(\Ga,q)=r(0,q)$
we get
$$\sum_{\al\in Q_+}q^{-d(\al,\nu)}P(\lM(\al,\nu),q)x^\al
=\sum_{\al\in\cN^I}a_{(\al,1)}(\Ga_*,q)x^\al=\frac{r(\nu,q)}{r(0,q)}.$$
\end{proof}

\begin{lmm}
We have
$$\sum_{\al\in\cN^I}r_{(\al,1)}(\Ga_*,q)x^\al=\frac {r(\nu,q)}{q-1}.$$
\end{lmm}
\begin{proof}
Let $T_*$ be the Tits form on $\cZ^{I_*}$ associated to the quiver
$\Ga_*$. Then, for any $\be\in\cZ^I$, we have
$T_*((\be,1))=T(\be)+1-(\nu,\be)$ and $T_*((\be,0))=T(\be)$. Any
partition of $(\al,1)$ is determined by a partition of
$\al$ as there exists just one partition of $1$. All these remarks
imply
\begin{multline*}\sum_{\al\in\cN^I}r_{(\al,1)}(\Ga_*,q)x^\al%
=\sum_{\tau\in\lP^I}\frac{q^{-1+(\nu,\tau_1)}}{1-q\inv}x^{|\tau|}%
\prod_{k\ge1}q^{-T(\tau_k)}\conjug{[\infty,\tau_k-\tau_{k+1}]}\\
=\sum_{\tau\in\lP^I}\frac{q^{(\nu,\tau_1)}}{q-1}x^{|\tau|}%
\prod_{k\ge1}{q^{-T(\tau_k)}} {\conjug{[\infty,\tau_k-\tau_{k+1}]}}
=\frac {r(\nu,q)}{q-1}.
\end{multline*}
\end{proof}

\begin{rmr}
Using the Hausel formula, the first Kac conjecture can be easily proved.
This was announced by Hausel in \cite{Haus1}. We give a brief idea of the proof.
The first Kac conjecture says that, for any $\al\in\cQ_+$, we have
$\dim\g_\al=a_\al(0)$, where \g is the Kac-Moody algebra associated to the
quiver \Ga and $a_\al(q)=a_\al(\Ga,q)$. Equivalently, for any
$\al\in Q_+$, we should have $\dim(U\g)_\al=m_{\al}(0)$, where
$\sum m_\al x^\al=\Exp(\sum a_\al x^\al)$. It holds
$\dim(U\g)_\al=\dim L(\nu)_{\nu-\al}$ for $\nu\gg0$ (this means that all
coordinates of $\nu$ are large enough).
From the results of Nakajima \cite{Nak2}
it follows that
\begin{multline*}
\dim L(\nu)_{\nu-\al}=\dim H^{\rm mid}_c(\lM(\al,\nu),\cC)\\
=q^{-d(\al,\nu)}P(\lM(\al,\nu),q)|_{q=0}
=\Big(\frac{r(\nu,0)}{r(0,0)}\Big)_\al=
(r(\nu,0)m(0))_\al.
\end{multline*}
 To show that $(r(\nu,0)m(0))_\al=m_\al(0)$ for $\nu\gg0$
one has to prove that $r(\nu,0)_0=1$ and
$r(\nu,0)_\be=0$ for $0<\be\le\al$ and $\nu\gg0$, which follows
from equation \eqref{eq:r_nu}.
\end{rmr}

%% file: section06.tex
\section{Combinatorics of fermionic forms}\label{sec:main}
All the generating functions of this section will be considered as
elements of the ring $R=\cQ(q)\lpoly{x_i,y_i| i\in I}$.
Given a function $f\in R$, we will denote by $f_\al$ its coefficient
by $y^\al$, where $\al\in\cN^I$. For any $\nu\in P$, we
define a ring homomorphism $S_\nu:R\ar R$ by
$$\sum f_\al y^\al\mto\sum q^{(\nu,\al)}f_\al y^\al.$$
Note that $S_\mu S_\nu=S_{\mu+\nu}$.
For any $\nu\in P$, define
$$p(\nu):=\sum_{\al\in Q_+}[\nu,\al]y^\al.$$
Define
$$p:=\sum_{\al\in Q_+}[\infty,\al]y^\al.$$

\begin{lmm}\label{lmm:equation:for:binomials}
We have $p(\nu)=p\cdot S_\nu\ub p$.
\end{lmm}
\begin{proof}
We can suppose that $\#I=1$.
Then, for any $n\in\cZ$,
$$p(n)=\sum_{k\ge0}[n,k]y^k
=\Exp\left(\frac{1-q^{n+1}}{1-q}y\right)$$
and
$$p=\sum_{k\ge0}[\infty,k]y^k=\Exp\left(\frac{y}{1-q}\right).$$
This implies
$$p\cdot S_n\ub p=\Exp\left(\frac y{1-q}+\frac{q^ny}{1-q\inv}\right)
=\Exp\left(\frac{1-q^{n+1}}{1-q}y\right)=p(n).$$

\end{proof}

Define
$$s(\nu):=\sum_{\tau\in\lP^I} x^{|\tau|}y^{\tau_1}\prod_{k\ge1}q^{T(\tau_k)}
\big[\nu-\sum_{l=1}^k \tau_l,\tau_k-\tau_{k+1}\big],$$
$$s:=\sum_{\tau\in\lP^I} x^{|\tau|}y^{\tau_1}\prod_{k\ge1}q^{T(\tau_k)}
\big[\infty,\tau_k-\tau_{k+1}\big].$$%

The following result is an analog of the Kleber algorithm
\cite{Kleber1}

\begin{prp}\label{prp:recursion formula}
For any $\nu\in P$ and $\al\in Q_+$, we have
\begin{enumerate}
    \item $s(\nu+\al)_\al=q^{T(\al)}x^\al ( p(\nu)\cdot
    s(\nu))_\al$.
    \item $s_\al=q^{T(\al)}x^\al ( p\cdot s)_\al$.
\end{enumerate}
\end{prp}
\begin{proof}
We will only prove the first formula. For any $\al\in Q_+$, there
is a bijection between the sets $\sets{\tau\in\lP^I}{\tau_1=\al}$
and
$$\bigcup_{0\le\be\le\al}\sets{\tau\in\lP^I}{\tau_1=\be},$$
where an element $\tau=(\tau_1,\tau_2,\dots)$ from the first set
is sent to the element $(\tau_2,\tau_3,\dots)$ from the second
set. This implies
\def\tempind#1{{\substack{\tau\in\lP^I\\ \tau_1=#1}}}
\enlargethispage{2\baselineskip}
\begin{multline*}
s(\nu)_\al=\sum_{\tempind\al}x^{|\tau|} \prod_{k\ge1}q^{T(\tau_k)}
\big[\nu-\sum_{l=1}^k
\tau_l,\tau_k-\tau_{k+1}\big]\\
=\sum_{0\le\be\le\al}\sum_{\tempind\be}x^{\al+|\tau|}
q^{T(\al)}[\nu-\al,\al-\be] \prod_{k\ge1}q^{T(\tau_k)}
\big[\nu-\al-\sum_{l=1}^k
\tau_l,\tau_k-\tau_{k+1}\big]\\
=q^{T(\al)}x^\al
\sum_{0\le\be\le\al}[\nu-\al,\al-\be]\sum_{\tempind\be}
x^{|\tau|}\prod_{k\ge1}q^{T(\tau_k)}\big[\nu-\al-\sum_{l=1}^k
\tau_l,\tau_k-\tau_{k+1}\big]\\
=q^{T(\al)}x^\al
\sum_{0\le\be\le\al}[\nu-\al,\al-\be]s(\nu-\al)_\be\\
=q^{T(\al)}x^\al\sum_{0\le\be\le\al}p(\nu-\al)_{\al-\be}s(\nu-\al)_\be
=q^{T(\al)}x^\al(p(\nu-\al)\cdot s(\nu-\al))_\al.
\end{multline*}
\end{proof}

\begin{rmr}
Note that the above formulas allow us to find the coefficients of
the term $x^\al y^\be$ in $s(\nu)$ or in $s$ recursively. On each
step either $\be\ne0$ and then \al decreases or $\be=0$ and then
we use the initial values,  that is the coefficients of the term $y^0$ in
$s(\nu)$ (respectively, in $s$), which equal $1$. This implies, in
particular, that the first formula uniquely determines the
functions $s(\nu)$ (with initial values as above).
\end{rmr}

The following theorem is our main combinatorial result

\begin{thr}\label{thr:maincombi}
We have $s(\nu)=s\cdot S_\nu\conjug s$.
\end{thr}
\begin{proof}
In view of the previous proposition and the remark after it, we
just have to prove that the functions $g(\nu):=s\cdot S_\nu\conjug
s$, $\nu\in\cZ^I$ satisfy the first formula of Proposition
\ref{prp:recursion formula}, \ie %
\begin{equation}\label{eq1}
q^{-T(\al)}(s\cdot S_{\nu+\al}\conjug s)_\al=x^\al( p(\nu)\cdot
s\cdot S_\nu\conjug s)_\al.
\end{equation}
The left hand side of this formula can be written in the form
\begin{multline}\label{eq2}
q^{-T(\al)}(s\cdot S_\al S_\nu\conjug s)_\al%
=\sum_{\be\le \al}q^{(\al,\be)-T(\al)}s_{\al-\be}\cdot (S_\nu\conjug
s)_\be\\
=\sum_{\be\le \al}q^{T(\be)-T(\al-\be)} s_{\al-\be}\cdot
(S_\nu\conjug s)_\be .
\end{multline}
Consider operators on the ring $R$
$$T:\sum f_\al y^\al\mto\sum q^{T(\al)}f_\al y^\al,\qquad
T\inv:\sum f_\al y^\al\mto\sum q^{-T(\al)}f_\al y^\al$$%
and the ring homomorphism $X:R\ar R$, $\sum f_\al y^\al\mto\sum
x^\al f_\al y^\al.$ Then the formula \eqref{eq1} in view of
\eqref{eq2} has the form
$$\sum_{\be\le\al}(T\inv s)_{\al-\be}(TS_\nu\conjug s)_{\be}
=X( p(\nu)\cdot s\cdot S_\nu\conjug s)_\al,$$%
\ie
\begin{equation}\label{eq3}
T\inv s\cdot TS_\nu\conjug s=X( p(\nu)\cdot s\cdot S_\nu\conjug s).
\end{equation}
It follows from Proposition \ref{prp:recursion formula} that
$s=TX( p\cdot s)$ or, equivalently,
$$T\inv s/Xs=X p.$$
This implies
$$T\conjug s/X\conjug s=\conjug{T\inv s/Xs}=\conjug {X p}=X\conjug p.$$
Now the formula \eqref{eq3} can be rewritten in the form
$$X( p\cdot S_\nu\conjug p)=X p(\nu).$$
But this follows from Lemma \ref{lmm:equation:for:binomials}.
\end{proof}

\begin{proof}[Proof of Theorem \ref{intro:thr2}]
For any $\nu\in P$, there is a ring homomorphism
$$\Phi_\nu:R'=\cQ(q)[y_i|i\in I]\lpoly{x_i|i\in I}\ar \cQ(q)\lpoly{x_i|i\in
I}$$ given by
$$\sum f_{\al,\be} x^\al y^\be\mto \sum q^{-(\nu,\be)}f_{\al,\be} x^{\al}.$$
Note that $s(\nu)$ and $s$ are in $R'$ and
$$n(\nu)=\Phi_\nu(s(\nu)),\qquad r(\nu)=\conjug{\Phi_\nu(s)}.$$
Hence
$$n(\nu)=\Phi_\nu(s\cdot S_\nu\conjug s)%
=\Phi_\nu(s)\cdot\Phi_0(\conjug
s)=\Phi_\nu(s)\cdot\conjug{\Phi_0(s)}=\conjug{r(\nu)}\cdot r(0).$$
This implies $n(\nu,q\inv)=r(\nu,q)r(0,q\inv)$. In particular, for
$\nu\in P_+$, we obtain
$$\sum_{\al\in Q_+}q^{-d(\al,\nu)}P(\lM(\al,\nu),q)x^\al
=\frac{r(\nu,q)}{r(0,q)}=\frac{n(\nu,q\inv)}{r(0,q)r(0,q\inv)}.$$
\end{proof}

%% file: section07.tex
\section{Fermionic Lusztig conjecture}\label{sec:ferm:Lusztig:conj}
This section is devoted to the proof of the fermionic Lusztig
conjecture. Recall that, for any quiver $(\Ga,I)$, we have defined
in Section \ref{sec:hausel formula} the polynomials
$a_\al(q)=a_\al(\Ga,q)$ and $m_\al(q)=m_\al(\Ga,q)$
($\al\in\cN^I$) counting the absolutely stable representation
(respectively, all representations) of \Ga over finite fields.

\begin{proof}[Proof of Theorem \ref{intro:thr1}]
We know that if \Ga is of type ADE then all the roots of $\g(\Ga)$
are real and therefore
$$a(q)=\sum_{\al\in Q_+} a_\al(q)x^\al=\sum_{\al\in\De_+}x^\al$$
does not depend on $q$, where $\De_+$ is the set of positive roots
of $\g(\Ga)$. This implies
$$r(0,q)r(0,q\inv)=\Exp\Big(\frac {a}{q-1}\Big)\Exp\Big(\frac {a}{q\inv-1}\Big)=\Exp(-a)
=\frac 1{m},$$%
where $m=\Exp(a)=\prod_{\al\in\De_+}(1-x^\al)\inv$. For $\al\in
Q_+$, we will identify $e^{-\al}$ with $x^\al$. It follows from
Theorem \ref{intro:thr2} that
$$\sum_{\al\in Q_+}q^{-d(\al,\nu)}P(\lM(\al,\nu),q)e^{-\al}%
=\frac{n(\nu,q\inv)}{r(0,q)r(0,q\inv)}=m\sum_{\la\in
P}n(\nu,\la,q)e^{\la-\nu}.$$%
Note that $m=\sum_{\al\in
\De_+}(1-x^\al)\inv=\prod_{\al\in\De_+}(1-e^{-\al})\inv$ and
therefore, for any $\la\in P$, we have $\ch M(\la)=e^\la m$ (see,
e.g., \cite[9.7.2]{Kac3}). This implies
$$m\sum_{\la\in P}n(\nu,\la,q)e^{\la-\nu}=e^{-\nu}\sum_{\la\in P}
n(\nu,\la,q)\ch M(\la)$$ and the theorem is proved.
\end{proof}